\documentclass[11pt]{article}
\usepackage{amsfonts,amssymb, graphicx,a4,bm,bbm}
\usepackage{color}

\newtheorem{theo}{Theorem}[section]
\newtheorem{lem}[theo]{Lemma}

\newtheorem{defi}[theo]{Definition}

\newtheorem{rema}[theo]{Remark}

\definecolor{Red}{cmyk}{0,1,1,0}

\let\a=\alpha \let\b=\beta  \let\d=\delta 
 \let\g=\gamma      
  \let\o=\omega     
\let\r=\rho  \let\t=\tau 
  
\let\D=\Delta \let\F=\Phi

\let\\=\noindent
\def\ee{\end{equation}}
\def\be{\begin{equation}}

\let\a=\alpha \let\b=\beta  \let\d=\delta 
 \let\g=\gamma      
  \let\o=\omega     
\let\r=\rho  \let\t=\tau 
  
\let\D=\Delta \let\F=\Phi

\let\\=\noindent
\def\ee{\end{equation}}
\def\be{\begin{equation}}

\def\0{\emptyset}

\title{A new bound on  the acyclic edge chromatic index}
\author{
\\
Paula M. S. Fialho, Bernardo N. B. de Lima$^*$, Aldo Procacci\\
\\
\small{ Departamento de Matem\'atica UFMG}
\small{ 30161-970 - Belo Horizonte - MG
Brazil}}

\begin{document}

\maketitle

\begin{abstract}
In this note we obtain a new bound for the  acyclic edge chromatic number $a'(G)$ of a graph $G$ with maximum degree $\D$  proving  that
$a'(G)\le 3.569(\D-1)$. To get this result we revisit  and slightly modify the method described
in  [Giotis,  Kirousis, Psaromiligkos and  Thilikos,  Theoretical Computer Science,  66: 40-50, 2017].
\end{abstract}

\vskip.2cm
{\footnotesize

}
\vskip.5cm

\vskip.5cm
\def\EA{{E_{\cal A}}}
\def\red{\color{Red}}

\section{Introduction}\def\GG{{\cal G}}
Let $G=(V,E)$ be a graph with vertex set $V$, edge set $E$ and maximum degree $\D > 1$, such that $|E|=m$. Given  $N\in \mathbb{N}$, let us denote $[N]=\{1,2,\dots,N\}$. A coloring of the edges of $G$ is a function $c: E \to [N]$. An edge coloring of $G$ is called  {\it proper} if no two adjacent edges receive the same color, and a proper edge coloring of $G$ is \textit{acyclic} if any cycle is colored with at least three colors. The minimum number of colors required such that a graph $G$ has at least one acyclic proper edge coloring is called the {\it acyclic edge chromatic number} of $G$  and is denoted by $a'(G)$. Given
a graph $G$ with maximum degree $\D$,  by Vizing Theorem a trivial lower bound for $a'(G)$ is $\D+1$.
The chronicle of the upper bound for the edge chromatic index  $a'(G)$ of a graph $G$ with maximum degree $\D$  goes back to the paper \cite{AMR} where the authors proved using the Lov\'asz local lemma (LLL) that there exists a constant $C\le 64$ such that  $a'(G)\le C\D$.
Since then efforts have  been done to lower the constant $C$. Molloy and Reed  showed  that $C\le 16$  in \cite{MR} using again Lov\'asz local lemma.
Fiam\v{c}ik \cite{Fi}, and later Alon et al. \cite{ASZ} have conjectured   that $a'(G)\le\D+2$. In \cite{ASZ}
 this conjecture is proved for graphs with girth  $g\ge 2000\D\ln\D$.
The upper bound on $a'(G)$ for general graphs with maximum degree $\D$ obtained by Molloy and Reed in 1998 was improved in 2012  by Ndreca et al. \cite{NPS}
who showed that  $a'(G)\le 9.62\D$ using an improved version of the Lov\'asz local lemma by Bissacot. et al. \cite{BFPS}. Only one year later
Esperet and Parreau \cite{EP} further improved this bound sensibly showing  that  $a'(G)\le 4\D$ by using following  the crucial observation.

\begin{lem}[Esperet-Parreau]\label{ep}
It is  possible to color  greedily  the edges  of a graph $G$ with maximum degree $\D$ using $N> 2(\D-1)$ colors
in such a way that the resulting coloring is proper and free of bichromatic cycles of length 4.
\end{lem}

\\In \cite{EP} authors manage to fit Lemma \ref{ep} into the so-called entropy compression method.
The entropy compression method (ECM) rests its basis on a sequential  algorithm inspired by the celebrated algorithmic version of Lov\'asz local lemma obtained by Moser and Tardos in 2010 \cite{Mo,MT} and it  can be applied to a wide class of graph coloring problems. Indeed, the ECM
has  been  successfully used to improve upper bounds for several  chromatic indices of bounded degree graphs previously obtained via the Lov\'asz local lemma (see e.g. \cite{GMP}, \cite{Pr}, \cite{DJKW}).

\\These achievements
have instilled the belief that the ECM is a tool  more efficient than the  Lov\'asz local lemma, even  its improved  version by Bissacot et al. \cite{BFPS}, as far as
graph coloring problems are concerned. Very recently however,  a further improvement on the notable Esperet-Parreau bound  for $a'(G)$
has been obtained by
Giotis et al. \cite{GKPT} who
show that $a'(G)\le 3.79(\D-1)$. Remarkably, the authors do not make use of the ECM. Rather they basically manage to accommodate Lemma \ref{ep} in   the standard Moser-Tardos scheme  (which in general is expected to give bounds identical to those obtained via the LLL).  This result suggests on one hand that
the  strong improvement obtained by Esperet and Parreau for $a'(G)$ is more due to Lemma \ref{ep} than to ECM and on the other hand that the Giotis bound could be further improved
being able to include in their scheme, besides Lemma \ref{ep},  some features of the ECM.

\\In this note we show that this is indeed the case.  By revisiting and slightly modifying  the method described by Giotis et al. we  get a further improvement
for the upper bound of $a'(G)$ obtaining that $a'(G)\le 3.569(\D-1)$.

\def\A{{\cal A}}

\section{\textsc{Color-Algorithm}}
\def\EA{{E_{\cal A}}}\def\red{\color{Red}}

\\Given a graph $G=(V,E)$ with maximum degree $\D$ such that $|E|=m$,
let $N= \left\lceil (2+ \gamma)(\D -1)\right\rceil+1$, where $\gamma$ is a positive number to be determinate later.
A pair of adjacent edges in $G$ is called a {\it cherry}. A $k$-cycle of $G$  is a cycle containing  $k$ edges.
 A partial coloring $w$ of $E$ is a function $w:E\to [N]_0$ where $[N]_0=[N]\cup\{0\}$ and
when  $w(e)=0$ we say that $e$ is uncolored.
We denote by $Y_{2k}$  the collection of all the $2k$-cycles of $G$, and we set $Y = \cup_{k\geq 3}{Y_{2k}}$.
Hereafter we suppose that an order is chosen in the sets $E$,  $V$, and  $Y$.


\\Given any edge $e \in E$ and any partial coloring $w$ of $E$ such that $w(e)=0$, let $D(e,w)\subset [N]$ be the set of available  colors for  the edge $e$ in order to avoid monochromatic cherries or bichromatic 4-cycle.
By Lemma \ref{ep}, 
we have that $\left|D(e,w)\right| \geq \left\lceil \gamma (\D-1)\right\rceil+1$.

\\We now describe a procedure, called \textsc{Color-Algorithm} which  colors (and eventually recolors) sequentially the edges of $G$.
Each discrete time $t\in \mathbb{N}$ \textsc{Color-Algorithm} colors (or recolors) an
edge is called an {\it instant} and  we denote by $w_t$  the coloring (or partial coloring) of $E$ at instant $t$.
Given a  coloring  of $E$ we say that an edge $e$ is {\it badly colored}  if there exists a cycle $C\in Y$ such that $e\in C$ and
$C$ is (properly) bichromatic. Conversely,   if there is no cycle $C\in Y$ such that $e\in C$ and
$C$ is bichromatic, $e$ is said to be {\it well colored}.

\vskip.2cm
\noindent
{\textsc{Color-Algorithm}}.
\vskip.2cm
\begin{enumerate}
  \item Color all edges $e \in E$ sequentially following the pre-fixed order in the following way: at each instant $t$ choose uniformly at random a number $r \in \{1, \cdots, \left\lceil \gamma (\D-1)\right\rceil+1\}$ and assign to $e$ the $r$-th smallest color in $D(e,w_{t-1})$.
  \item  While there is a badly colored edge, let $e$ be the largest edge among them and let $C$ be the smallest bichromatic cycle such that $e$ is one of its edges, and do
  \item   \textsc{\textsc{Recolor}}$(e, C)$.
  \item End while.
  \item Output current evaluation.
\end{enumerate}

\vskip.5cm

\\{ \textsc{\textsc{Recolor}}$(e,C)$}

\begin{enumerate}
  \item Let $f_1$ be the edge of $C$ to receive the color
  that it has in this phase
  at the earliest instant among all the edges in $C$, and let $f_2$ be the edge in $C$ among those with opposite parity w.r.t $f_1$ to receive the color that it has in this phase earliest.
  Define $ S(C)= (f_1, f_2)$.
  \item  \textsc{Recolor} all the edges in $C\setminus S(C)$ sequentially (according  the pre-fixed order in $E$) assigning to the edge
   recolored at instant $t$ the $r$-th smallest color in $D(f,w_{t-1})$, where $r$ is chosen uniformly at random in $ \{1, \cdots,$ $ \left\lceil \gamma (\D-1)\right\rceil+1\}$.
  \item While there exists an edge in  $ C\setminus S(C)$ which is badly colored, let $e'$ be  the largest of these edges  and  let $C'$ be the smallest bichromatic cycle such that $e'$ is one of its edges, and do
  \item \textsc{\textsc{Recolor}}$(e',C')$.
  \item End while.
\end{enumerate}

\vskip.2cm
\\Note that for any $t<m$ $w_t$ is a partial proper coloring
without bichromatic 4-cycles while,
for any $t\ge m$, $w_t$ is a  proper coloring
without bichromatic 4-cycles.
\\A {\it step} of \textsc{Color-Algorithm}  is the procedure described in Line 2 of \textsc{Recolor}$(e,C)$.
\\A {\it phase } of \textsc{Color-Algorithm} is the collection of steps made by \textsc{Color-Algorithm}
during a  call of\textsc{ Recolor}$(e,C)$ in Line 3 of \textsc{Color-Algorithm}. The  {\it root} of the phase is its initial step.
 The set $S(C)$ defined in Line 1 of \textsc{Recolor}$(e,C)$ will be called a {\it seed} of $C$.
The  {\it record}
of the algorithm is the list   $\mathcal{L}=((e_1,C_1), (e_2,C_2),\dots, $   constituted by  steps done by the algorithm during its execution. According to the prescriptions described above, $\mathcal{L}$ is a random variable
determined by the random samplings performed by the algorithm in each  step. If  $\mathcal{L}$ is finite, i.e. if $|\mathcal{L}|=n$ for some $n\in \mathbb{N}$,  then the algorithm terminates having performed $n$ steps and  produces an acyclic edge coloring of $G$.
Let us define
\begin{equation}\label{p}
  P_n= {\mathbb P}(|\mathcal{L}|=n).
\end{equation}
In other words $P_n$ is the probability that \textsc{Color-Algorithm} runs $n$ steps.

\begin{rema}\label{r46}
Our procedure \textsc{Recolor}$(e,C)$  is similar but not identical to that described in \cite{GKPT}.In \cite{GKPT} at each step of  \textsc{Recolor}$(e,C)$  all edges in $C$ are recolored
while we recolor only  the edges in $C \setminus S(C)$.
\end{rema}

\\Let us now prove three key properties of \textsc{Color-Algorithm}.
\vskip.2cm
\begin{lem}\label{vinseed}
  In any call of \textsc{Recolor}$(e,C)$ the edge $e$ does not belong to $S(C)$ and thus  is always resampled.
\end{lem}
{\it Proof:} Consider first a root call of \textsc{Recolor}(e,C) in Line 3 of \textsc{Color-Algorithm}. Note that   the pair  $(e,C)$  chosen in Line 2 of \textsc{Color-Algorithm} is such that $e$ is the
largest edge in the cycle $C$  and therefore $e$   was colored after all the other edges of $C$, which means that $e$ never  belongs to the seed $S(C)$. Consider now the recursive call of \textsc{Recolor}$(e',C')$ in Line 4 of \textsc{Recolor}$(e,C)$. Observe that in this case the edges in $C \setminus S(C)$ have been resampled and $e'$ was taken as the largest edge in $C \setminus S(C)$ that is still in a bichromatic cycle, and $C'$ is the smallest bichromatic cycle such that $e'$ is one of its edges. If $C'=C$ it means that the same cycle is still bichromatic, so the seed now is the same as before and $e'$ is not in the seed. If $C' \neq C$,
then $e'$ is the largest among the egdes of $C' \cap (C \setminus S(C))$ which have been resampled (while the other edges in $C'$ have not).
In conclusion,  at the beginning of the step $(e',C')$,   $e'$ is the  edge resampled at the latest instant among those belonging to $C'$ and thus will not belong to $S(C')$.
$\Box$

\vskip.2cm

\vskip.2cm

\begin{lem}\label{Agood}
Consider any call of \textsc{Recolor}$(e,C)$, let $S(C)$ be the seed chosen at the beginning of the call and let $X$ be the set of all well colored edges at the beginning of this call. If the call \textsc{Recolor}$(e,C)$ ends, then all the edges in $X\cup \left( C \setminus S(C) \right)$ are well colored.
\end{lem}
{\it Proof:} According to the algorithm, if \textsc{Recolor}$(e,C)$ ends then all edges in $C \setminus S(C)$ are not in a bichromatic cycle. So we just need to prove that no edge of $X$ is in a bichromatic cycle at the end of \textsc{Recolor}$(e,C)$. By contradiction, assume that \textsc{Recolor}$(e,C)$ lasts for  $n$ steps and  $f \in  X$ belongs to a bichromatic cycle $B$ after the last step. This cycle  $B$ was  not bichromatic at the beginning of \textsc{Recolor}$(e,C)$, it has some of its edges resampled during the execution of \textsc{Recolor}$(e,C)$ and it is bichromatic when  \textsc{Recolor}$(e,C)$ ends. Therefore  there must exists a last step $s\le n$ of \textsc{Recolor}$(e,C)$  such that $B$ is not bichromatic at step $s-1$, becomes bichromatic at step $s$ and stays bichromatic during the remaining $n-s$ steps of \textsc{Recolor}$(e,C)$. According to \textsc{Recolor}$(e,C)$, there must be  a cycle $B'$ and an edge $f'$
such that the process \textsc{Recolor}$(f',B')$ was called at step $s-1$ of \textsc{Recolor}$(e,C)$ and $B$  shares  at least an edge $g$ with the set $B'\setminus S(B')$, and $B$ becomes bichromatic as soon as the edges in $B' \setminus S(B')$ were recolored. However, the algorithm says that the edge $g \in \{B'\setminus S(B') \}\cap B$ must not be in a bichromatic cycle at the end of the call of \textsc{Recolor}$(f',B')$, therefore at  some step $s'>s$, $B$, which contains $g$, must cease to be bichromatic. We have reached a contradiction.  $\Box$

\begin{lem}\label{mtrees}
\textsc{Color-Algorithm} performs at most $m=|E|$ phases.
\end{lem}
{\it Proof}. Consider two phases $l$ and $s$, with $l<s$, generated by an execution of \textsc{Color-Algorithm} and let $(e_l, C_l)$ and $(e_s, C_s)$ be the pairs resampled at their root steps respectively.  By Lemma \ref{Agood}, all edges in $(C_l)\setminus S(C_l)$ are not in a bichromatic cycle when phase $l$  ends (and at the beginning of any successive phase). In particular, since by Lemma \ref{vinseed} $e_l\in (C_l)\setminus S(C_l)$, $e_l$ is not in a bichromatic cycle and thus $e_l\notin C_s$.  In conclusion  $e_l \neq e_s$. $\Box$

\section{Witness forests}
We will  now associate to the record $\mathcal{L}$ of \textsc{Color-Algorithm} a labeled forest formed by plane rooted trees whose internal vertices are labeled with pairs  $(e,C)$ belonging to $\cal L$.

\\Suppose that the algorithm performs $r$ phases and that during the phase $s$ the algorithm performs $n_s$ steps in such a way that
the record of the algorithm is
\begin{equation}\label{forest}
 \mathcal{L}\; =\;\{(e^1_1,C_1^1), \dots, (e^1_{n_1},C_{n_1}^1),(e^2_1,C_1^2), \dots, (e^2_{n_2},C_{n_2}^2), \dots, (e^r_1,C_1^r), \dots, (e^r_{n_r},C_{n_r}^r)\}.
\end{equation}
We will now associate a labeled rooted tree $\t'_s$ to each phase $s$, $1 \leq s \leq r$. Let
\begin{equation}\label{t1}
  (e^s_1,C_1^s), \dots,({e^{s}_{i}},C^s_i), \cdots, (e^s_{n_s},C_{n_s}^s),
\end{equation}
be the pairs recolored at phase $s$. We construct the tree $\t'_s$ in the following way.

\vskip.2cm
\\a) the root of  $\t'_s$  has label  $(e^s_1,C_1^s)$.

\\b)  For $i=2,\dots, n_s$, we proceed  by checking  if $({e^{s}_{i}},C^s_i)$ is such that ${e^{s}_{i}} \in (C^{s}_{i-1} \setminus S(C^s_{i-1}))$,

-\\ if yes, we add $({e^{s}_{i}},C^s_i)$ as a child of $({e^{s}_{i-1}},C^s_{i-1})$, \\

- \\if no, we go back in (\ref{t1}) checking the {\it ancestors}  of the (already added)  vertex with label $({e^{s}_{i-1}},C^s_{i-1})$ until find a pair $({e^{s}_{j}},C^s_j)$, with $j<i$, such that ${e^{s}_{i}} \in C^s_{j} \setminus S(C^s_{j})$, and we add $({e^{s}_{i}},C^s_i)$ as a child of $({e^{s}_{j}},C^s_{j})$.

\\Note that $\t'_s$ has $n_s$ vertices
(leaves included) with labels $({e^{s}_{i}},C^s_i)$ with $i=1,\dots, n_s$.
Observe moreover that, by Lemma \ref{Agood}, the pair $(e^{s+1}_1,C_1^{s+1})$ is  such that $e^{s+1}_1\notin C^l_i\setminus S(C^l_i)$ for all $i\in [n_l]$ and for all $l\le s$ ,
so we build a new tree $\t'_{s+1}$ with root $(e^{s+1}_1,C_1^{s+1})$ following the same rule described to build $\t'_s$.

\\The  forest ${ F}'=\{\t'_1,\dots, \t'_r\}$ so obtained uniquely associated to the record $\mathcal{L}$ is such  that, for each $s\in [r]$,   $\t'_s$ is a rooted plane tree with  $n_s$  vertices and each vertex of $\t'_s$ has label $(e,C)$ where $e\in E$ and $C\in Y$. Note that, by Lemma \ref{mtrees} we have that $r\le m$ and thus the   forest ${ F}'$ contains at most $m$ trees.

\\Note also that in each tree $\t'_s$ of $F'$ the list of labels of the vertices of $\t'_s$,  ordered according to the depth-first search, coincides with  the list  (\ref{t1}).

\begin{rema}\label{inje}
By construction, the correspondence ${\cal L}\mapsto F'$ is an injection.
\end{rema}

Let us show three key properties of the trees in the forest  ${ F}'$.

\begin{lem}\label{dr}
Given a tree $\t\in { F}'$, the following holds:

\begin{enumerate}

\item[\rm a)] Let the vertex $u$ be a child of the vertex $v$ and let $(e,C)$ and $(f, B)$ be their labels respectively. Then $e \in B\setminus S(B)$.

\item[\rm b)] Let the vertices $v$ and $v'$ be the $i^{th}$ and the $j^{th}$  siblings in $\t$, with $i<j$ in the natural order of the vertices of $\t$ induced by the steps of the algorithm (i.e. the depth-first search order of $\t$), and let $(e_i,C_i)$ and $(e_j, C_j)$ be their labels  respectively, then $e_i\neq e_j$.
\item[\rm c)] Let $v$ be a vertex  and let $(e,C)$ be its label, then the vertex $v$ has at most $|C|-2$ children.
\end{enumerate}
\end{lem}

\\{\it Proof.}
Item a) is trivial by construction of the algorithm. Item c) follows trivially from item b). Let us thus prove item b). For $k$ such that $i\le k<j$, let $(e_k,C_k)$ be the label of the $k^{th}$ sibling. By Lemma \ref{Agood} when Recolor$(e_k,C_k)$ ends all the edges label $e_i,\dots, e_k$ are  not in a bichromatic cycle and also the edges in  $(C_i\setminus S(C_i))\cup\dots\cup (C_k\setminus S(C_k))$.
Therefore $e_j$ can not be in the set $\{e_i,e_{i+1},\dots, e_{j-1}\}$. $\Box$
\vskip.2cm

\\Let ${E}_{F'}=\{e^1_1,e^2_1, \dots,  e^r_1\}$. I.e., for $s\in [r]$, $e^s_1$ is the edge label of the root of the tree $\t'_s\in F'$.
By Lemma \ref{mtrees}, the edges  in  ${E}_{F'}$ are all distinct.

\begin{defi}[Witness forest]
  Given  the record of \textsc{Color-Algorithm} $\mathcal{L}$, the  witness forest $F$ associated to $\mathcal{L}$ is built  starting from $F'$ in the following way:

  \\1) Add to the forest $F'$ as many isolated vertices  as the edges which are in $ E \setminus { E}_{F'}$, and  give to these isolated vertices the label $(e, \emptyset)$ for all $e \in { E}\setminus { E}_{F'}$.

  \\2) For each vertex $v$ of the forest $F'$ with cycle label $C_v$ which has less than $|C_v|-2$ children do the following.  Let $H_v$ de the set of edges in $C_v\setminus S(C_v)$ which are not edge labels of the children of $v$. For each $e \in H_v$ add to $v$ a leaf with label $(e,\emptyset)$ in such a way that $v$ has now exactly  $|C_v|-2$ children.
\end{defi}

\\The  new labeled forest $F$ so obtained uniquely  associated to the random variable $\mathcal{L}$ by the prescriptions described above is called the {\it witness forest} produced by \textsc{Color-Algorithm}.  This witness forest $F$ (a random variable)  has the following properties.

\noindent
\vskip.2cm
{\bf Properties of the witness forest $F$}.
\vskip.2cm
\begin{enumerate}

\item $F$  is constituted by  exactly $|E|=m$ labeled rooted trees $\t_1,\dots,\t_m$ (some of which are just isolated vertices).

\item  Each internal vertex $v$ of $\t\in F$ carries a label $(e_v,C_v)$ where $e_v\in E$ and $C_v\in Y$ while each leaf $\ell$ of $\t$
carries a label $(e_\ell,\0)$ where $e_\ell\in E$.

\item Let $v$ be an internal vertex of $\tau \in F$ with label $(e_v, C_v)$ and let $u$ be a child of $v$. Then the edge label of $u$
belongs to $C_v$.

\item Let the vertices $v$ and $v'$ be  siblings in $\t\in F$  and let $(e_v,C_v)$ and $(e_{v'}, C_{v'})$ be their labels  respectively, then $e_v\neq e_{v'}$.
\item Let $v$ be an internal vertex of $\t\in F$   and let $(e_v,v)$ be its label, then the vertex $v$ has exactly $|C_v|-2$ children.
\end{enumerate}
Let  $\mathfrak{F}_n$ be the set of labeled   forests satisfying properties 1-5 with $n$ internal vertices  and let
$\mathfrak{F}=\cup_{n\ge 0} \mathfrak{F}_n$.

\\Recalling Remark \ref{inje}, it is important to stress that the map ${\cal L}\mapsto F$ is an injection. In other words,  distinct records ${\cal L}_1$ and ${\cal L}_2$ necessarily produce distinct witness forest $F_1$ and $F_2$.
\\Therefore, since \textsc{Color-Algorithm} lasts $n$ steps if and only if the witness forest associated to the record $\cal L$ of \textsc{Color-Algorithm} has $n$ internal vertices, the probability
$P_n$ defined in (\ref{p}) can be written as
\begin{equation}\label{wf}
P_n    = {\mathbb P}(\mbox{the witness forest produced by \textsc{Color-Algorithm} $\cal L$ has} \ n \ \mbox{internal vertices}).
\end{equation}

\section{Algorithm \textsc{ColorVal}}

\\We now introduce a validation algorithm  following closely Giotis et al. \cite{GKPT}.

\begin{defi}
A  pair of adjacent edges $(e,f)$ is ordered if $e=\{v,u\}$, with  $v<u$, and $f$ contains $u$. Given an ordered pair $(e,f)$ and an integer $k\ge 3$ we say that the triple $(e,f,k)$ is admissible in $G$ if there exists at least one cycle $C$ in $G$ containing the ordered pair $(e,f)$  such that $|C|=2k$.
\end{defi}

\def\SS{{\mathcal S}}
\begin{defi}
Given $n\in \mathbb{N}$,
  we say that $S=\{(e_1^1, e_1^2, k_1), \cdots, (e_n^1, e_n^2, k_n)\}$ is an admissible sequence of $G$ if, for all $i=1, \cdots, n$, the triple  $(e_i^1, e_i^2, k_i)$ is admissible in $G$. We denote by $\mathfrak{S}_n$ the set of all admissible sequences of $G$ constituted by $n$ triples and we set $\mathfrak{S}=\cup_{n\ge 0}\mathfrak{S}_n$.
\end{defi}

\begin{rema}\label{cru}
Given a proper edge coloring of $G$ and given a pair of adjacent edges $(e,f)$ of $G$, there exist at most one bichromatic cycle $C$ containing $e$ and $f$.
\end{rema}

\\We now describe an algorithm, called \textsc{ColorVal},  whose
 input is an admissible sequence $S=\{(e_1^1, e_1^2, k_1), \cdots,$ $ (e_n^1, e_n^2, k_n)\}$. \textsc{ColorVal}
 first  colors  sequentially all edges of $G$ and then recolors (always sequentially) some of the edges of $G$.
As before, each discrete time $t\in \mathbb{N}$ \textsc{ColorVal} colors (or recolors) an
edge is called an {\it instant} and  we denote by $w_t$  the coloring (or partial coloring) of $E$ at instant $t$.


\vskip.2cm
\noindent
\textsc{ColorVal(S)}.
\vskip.2cm
Given the admissible sequence $S=\{(e_1^1, e_1^2, k_1), \cdots,$ $ (e_n^1, e_n^2, k_n)\}$
\begin{enumerate}
  \item   Color all edges $e \in E$ sequentially following the pre-fixed order in the following way: at each instant $t\le m$ choose uniformly at random a number $r \in \{1, \cdots, \left\lceil \gamma (\D-1)\right\rceil+1\}$ and assign to $e$ the $r$-th smallest color in $D(e,w_{t-1})$.
  \item  For $i= 1, \cdots, n$, do
  \item  If there is a bichromatic $2k_i$-cycle containing $e_i^1$ and $ e_i^2$, let $C_{i}$ be the unique such cycle. If there is not a such bichromatic cycle, let $C_i$ be an arbitrary (e.g. the smallest) cycle containing $e_i^1$ and $ e_i^2$  such that $|C_i|=2k_i$. \textsc{Recolor} all the edges $ e \in C_{i} \setminus S({C_i})$ following the pre-fixed order and choosing at each instant $t>m$ a color uniformly at random in $D(e,w_{t-1})$.
  \item End for.
\end{enumerate}

\\The procedure described at line 3 of ColorVal is called a {\it step} and  $S({C_i})$ is the seed of ${C_i}$ defined as in Line 1 of \textsc{Recolor}$(e,C)$. Of course,
if $S=\{(e_1^1, e_1^2, k_1), \cdots, (e_n^1, e_n^2, k_n)\}$ is the input for \textsc{ColorVal},  its execution will perform exactly $n$ steps.
Moreover, with such an input, \textsc{ColorVal} assigns a color to an edge exactly $m+\sum_{i=1}^n(2k_i-2)$ times, i.e., its execution will last for  exactly $m+\sum_{i=1}^n(2k_i-2)$ instants.

\begin{rema}\label{r46}
Our algorithm \textsc{ColorVal} is similar but not identical to that described in \cite{GKPT}. The difference
comes from the fact that  in \cite{GKPT} at each step $i$  of \textsc{ColorVal}  all edges in $C_i$ are recolored
while our \textsc{ColorVal} recolors only  the edges in $C_i \setminus S(C_i)$.
\end{rema}

\\Note that the algorithm \textsc{ColorVal} produces as output an unique cycle sequence $\mathcal{C}= \{C_1, \cdots, C_n$\}. We say that the algorithm \textsc{ColorVal} is {\it successful} if all cycles in $\mathcal{C}$ were chosen as bichromatic. The lemma below furnishes  an upper bound
for the probability that \textsc{ColorVal} is successful.
\begin{lem}\label{l5}
  Given an admissible sequence $S= \{(e_1^1, e_1^2, k_1), \cdots, (e_n^1, e_n^2, k_n)\}$, it holds that
  \begin{equation} \label{cova}
    \mathbb{P}(\mbox{\textsc{ColorVal} is successful in S}) \leq \left(\frac{1}{\lceil\gamma(\Delta -1)\rceil+1}\right)^n \prod_{s=1}^{n}{\left(1-\left(1 - \frac{1}{\lceil\gamma(\Delta -1)\rceil+1} \right)^{\Delta-1}\right)^{2k_s -3}}.
  \end{equation}
\end{lem}
\\The proof of this lemma  is  identical  to the proof of Lemma 5  in \cite{GKPT}, where Remark \ref{cru} is crucial. This is so despite the fact that
our \textsc{ColorVal} is slightly different from that described in \cite{GKPT} (recall Remark \ref{r46}). The key point is that, given a cycle
$C$ belonging to the output $\cal C$ of \textsc{ColorVal}, the probability  that  the edges of the seed $S(C)$ (which Giotis et al. call the {\it early edges of $C$}) receive the color that makes $C$ bichromatic is (by definition) equal to one
and thus they  do not play any role in the estimate (\ref{cova}).

\begin{rema}\label{lemma5b}
\\Using that $\lceil\gamma(\Delta -1)\rceil+1\ge \gamma(\Delta -1)+1$ and the inequality $ 1 - \frac{1}{1+x} > e^{-\frac{ 1}{x}} $, valid for all $x > 0$,  the bound ( \ref{cova})can be simplified as follows.
\be
\mathbb{P}(\mbox{\textsc{ColorVal} is successful in S}) \le  \frac{1}{ (\D -1)^n} {\prod_{s=1}^{n}{\left(\frac{1}{\gamma}\left( 1- e^{\frac{-1}{\gamma}}\right)^{2k_s -3}\right)}}.
\ee
\end{rema}

\vskip.2cm

\def\FF{{\mathcal F}}

\\Given $n\in \mathbb{N}$, let $F\in \mathfrak{F}_n$ with internal vertices (ordered according to the depth-first search) with labels $(e_1,C_1),\dots,(e_n,C_n)$.
 We associate uniquely to $F$  the admissible sequence $S= \{(e_1^1, e_1^2, k_1), \cdots, (e_n^1, e_n^2, k_n)\}$ obtained by letting $e_i^1$ be $e_i$, $e_i^2$ be the neighbor of $e_i^1$ in $C_i$ such that $(e^1_i,e_i^2)$ is ordered, and $k_i = \frac{|C_i|}{2}$, for all $i=1, \cdots, n$. Given $S\in \mathfrak{S}$, let ${\FF}_S$ be the set of all witness forests in $\mathfrak{F}$ such that $S$ is the corresponding admissible sequence associated.

\begin{lem}\label{lem}
  Given an  admissible sequence $S\in \mathfrak{S}$, we have that
  \begin{equation}\label{colorval}
     \sum_{F \in \mathcal{F}_S}{\mathbb{P}(\mbox{\textsc{Color-Algorithm} produces the witness forest} \ F)} \le \mathbb{P}(\mbox{\textsc{ColorVal} is sucessful in S}).
  \end{equation}
\end{lem}

\\{\it Proof:} For any $F \in \FF_S$, let $Z_F$ be the event ``\textsc{Color-Algorithm} produces the witness forest $F$ ". Observe that
\begin{equation}
  \mathbb{P}\left(\bigcup_{F\in \FF_S}{Z_F}\right)= \sum_{F \in \FF_S}{\mathbb{P}(\mbox{\textsc{Color-Algorithm} produces the witness forest} \ F)},
\end{equation}
since the events $Z_F$ are mutually disjoint. On the other hand,  $\mathbb{P}\left(\bigcup_{F\in \FF_S}{Z_F}\right)$ is the probability that at least one forest $F \in \FF_s$ is the witness forest produced by \textsc{Color-Algorithm}. Note that if all the random choices made by an execution of \textsc{Color-Algorithm} that produces a record $\cal L$ such that its associated witness forest $F$ is in $\FF_S$  are also made by the algorithm\textsc{ColorVal} with input $S$, then \textsc{ColorVal} is successful. So,
\begin{equation}
\mathbb{P}\left(\bigcup_{F\in \FF_S}{Z_F}\right) \le \mathbb{P}(\mbox{ColorVal is sucessful in S}).
\end{equation}
$\Box$
\vskip.2cm

\section{A new bound for acyclic coloring}
This last section is devoted to the proof of the following theorem.
\begin{theo}\label{end}
The acyclic edge chromatic index $a'(G)$ of a graph $G$ with maximum degree $\D$ admits the following bound
\be\label{thm}
a'(G)\le 3.569(\D-1).
\ee
\end{theo}

\\The strategy to prove Theorem \ref{end} will be to show that
the probability (\ref{wf}) that \textsc{Color-Algorithm} lasts for $n$ steps decays exponentially with $n$, which implies that
 \textsc{Color-Algorithm} terminate almost surely returning an acyclic edge coloring. If \textsc{Color-Algorithm} lasts for $n$ steps  then it
 produces a witness forest  with $n$ internal nodes. Recall  that, if the internal vertex $v$  of the witness forest has cycle label $C_v$,  then this vertex has $|C_v|-2$ children.  Given the record $\cal L$ of \textsc{Color-Algorithm} such that $|{\cal L}|=n$ and given the witness forest $F\in \mathfrak{F}_n$ associated to $\cal L$ (i.e with $n$ internal vertices), by removing all its labels we obtain an unlabeled witness forest $\Phi$. We call $\Phi$ the
{\it associated unlabeled witness forest} produced by \textsc{Color-Algorithm}.
 This unlabeled forest $\Phi$  is constituted by $|E|=m$ plane rooted trees,  it has in total $n$ internal vertices and it is such that
each internal unlabel  vertex $v$ of $\Phi$ has $|C_v|-2$  children  with $C_v$ being the label of the corresponding vertex of $F$.
Let  ${\cal F}_n$ be the set of unlabeled forests constituted by $|E|=m$ trees with $n$ internal vertices and such that each internat  vertex $v$ has
a number of children in the set   $\{2k-2\}_{k\ge 3}$. For $\F \in {\cal F}_n$ let us define
\begin{equation}
  P_\F = {\mathbb P}( \F \ \mbox{is the associated unlabeled witness forest produced by \textsc{Color-Algorithm}}).
\end{equation}

\begin{lem}
Let $\Phi\in {\cal F}_n$ with internal vertices $v_1,\dots, v_n$. Let, for $s=1,2,\dots, n$,  $2k_s-2$ be the number of children of the internal vertex $v_s\in \Phi$. Then
  \begin{equation}\label{fi}
    P_\F \leq {\prod_{s=1}^{n}{\left(\frac{1}{\gamma}\left( 1- e^{\frac{-1}{\gamma}}\right)^{2k_s -3}\right)}}.
  \end{equation}
\end{lem}

\\{\it Proof:} Given $\F\in {\cal F}_n$, $\F$ will be the unlabeled witness forest of an execution of \textsc{Color-Algorithm} if and only if this execution produces a record $\mathcal{L}= (e_1, C_1), \cdots, (e_n, C_n)$ that can be associated to $\F$. Then, let us check what are the possibilities for $\mathcal{L}$ and calculate the probability that one execution produces a such sequence $\mathcal{L}$.

\\First of all observe that, given $\F\in {\cal F}_n$, the edge-label  $e_1$  of  the first pair of the sequence $\mathcal{L}$
is uniquely determined. Indeed, if the last  $i$ trees of $\Phi$ are isolated roots ({recall that the edges-label of all vertices  are selected as the largest edge in a bichromatic cycle}), then  $e_1=e_{m-i}$   is the sole possible edge-label of the root $\r_1$ of the last non trivial tree of $\Phi$. Now we need that a cycle $C_1$,  such that $e_{1}\in C_1$,  is bichromatic. The unlabeled forest $\Phi$ tells us
that such $C_1$ must have $2k_1=s_{\r_1}+2$ edges, where $s_{\r_1}$ is the number of  children of  the root $\r_1$ of the last non trivial tree of $\Phi$.
Let us choose one edge  $e_1^2$ among those incident to the largest vertex of $e_1$; we have $(\D-1)$ possibilities for $e_1^2$. Now we have the triple $(e_1, e_1^2, k_1)$ and $C_1$ is one
of the cycles containing  $e_1$ and $e_1^2$ and of size $2k_1$.

\\For each possibility of $C_1$, the  two edges of $C_1$ with opposite parity receiving  their colors earliest form the seed $S(C_1)$, and thus we know what are the edges label for all the children of $(e_1, C_1)$ in $\F$.
The next edge label $e_2$ in $\mathcal{L}$ is chosen as the label of the last child of $(e_1, C_1)$ that is not a leaf. Again we need  a cycle $C_2$ containing $e_2$ to be bichromatic, such that $C_2$ must have $2k_2=s_{\r_2}+2$ edges, where $s_{\r_2}$ is the number of children that this first child of $(e_,C_1)$ has. To determine  this cycle,  let us choose an  edge $e_2^2$ incident to the largest vertex of  $e_2$; we have at most $(\D-1)$ possibilities. Observe then that now we have an admissible triple $(e_2, e_2^2, k_2)$ such that we need that there exists a bichromatic cycle $C_2$ such that $e_2$ and $e_2^2$ are ordered neighbors and $|C_2|=2k_2$.

\\We continue in this way, following the structure of $\F$: when a leaf is reached, we go back to the last internal node (in a depth-first search way), and look to the next child that is not a leaf; the edge that labels this child will be the next edge-label in $\mathcal{L}$. When a tree of the forest $\F$ is exhausted, we go to the next root. As $\F$ has $n$ internal nodes, then we will have a factor $(\D-1)^n$ (the possibilities for each $e_i^2$).

\\Observe that at each step we will have a triple $(e_s, e_s^2,k_s)$ such that we want that there exists a bichromatic cycle $C_s$ such that $|C_s|=2k_s$ and $e_s$ and $e_s^2$ are ordered neighbors. Then, in fact what we need is that \textsc{ColorVal} is successful with the entry $S=\{(e_s, e_s^2,k_s)\}_{s=1}^{n}$.

\\By Lemma \ref{lem} we have that the probability that at least one forest $F \in F_S$ becomes the witness forest for \textsc{Color-Algorithm} is bounded by the probability that \textsc{ColorVal} is successful in $S$. Recalling Remark \ref{lemma5b},

\begin{equation}
  \mathbb{P}(\mbox{ColorVal is sucessful in S}) \le \frac{1}{ (\D -1)^n} {\prod_{s=1}^{n}{\left(\frac{1}{\gamma}\left( 1- e^{\frac{-1}{\gamma}}\right)^{2k_s -3}\right)}},
\end{equation}
counting the factor $(\D-1)^n$, we have that
\begin{equation}
  P_\F \leq {\prod_{s=1}^{n}{\left(\frac{1}{\gamma}\left( 1- e^{\frac{-1}{\gamma}}\right)^{2k_s -3}\right)}}.
\end{equation}
$\Box$
\vskip.2cm
\\We have clearly that the probability $P_n$ (see (\ref{wf})) that \textsc{Color-Algorithm} lasts $n$ steps is bounded by
$$
P_n\le  \sum_{\F\in {\cal F}_n}{P_\F}.
$$
\\To estimate  $\sum_{\F\in {\cal F}_n}{P_\F}$, observe that every forest $\F\in {\cal F}_n$ is constituted by $m$ trees $\t_1,\dots,\t_m$
with $n_1,\dots, n_m$ internal vertices such that $n_i\ge 0$ for all $i=1,\dots, m$ and such that $n_1+n_2+\dots+n_m=n$.
Note  the number of children of the internal vertices  of any  $\t\in \Phi$ takes values in the set $I=\{4,6,8,\dots\}$. Let us denote by $\cal T$ the set of plane trees with
number of children of the internal vertices taking  values in the set $I$ and let ${\cal T}_n$ be the set of  the tree in  $\cal T$ with $n$ internal vertices.

\\Let us denote  shortly,  for $k\in \{3,4,\dots\}$,
\be\label{wk}
w_k={\left(\frac{1}{\gamma}\left( 1- e^{\frac{-1}{\gamma}}\right)^{2k -3}\right)}.
\ee
For  a tree $\t\in \cal T$, let $V_{\t}$ be the set of its internal vertices. If $v\in V_{\t}$, we denote by $s_v$ the number  of its children and let $k_v={s_v+2\over 2}\in \{3,4,\dots\}$.   Then  define the weight of $\t$ as
$$
\o(\t)= \prod_{v\in V_\t}w_{k_v},
$$
and, for a given $n\in \mathbb{N}$, let
$$
Q_n=\sum_{\t\in {\cal T}_n}\o(\t).
$$
Therefore, the probability that  \textsc{Color-Algorithm} lasts $n$ steps is bounded by
\be\label{penne}
P_n\le  \sum_{n_1+\dots +n_m=n\atop n_i\ge 0}Q_{n_1}\dots Q_{n_m}.
\ee
It is now easy to check that  $Q_n$ is defined by the recurrence relation
\be\label{recur}
Q_n=\sum_{k \ge 3}w_k\sum_{n_{1}+\dots +n_{2k-2}=n-1\atop n_1\ge 0, \dots ,n_{2k-2}\ge 0}
Q_{n_1}\dots Q_{n_{2k-2}},
\ee
with $Q_0=1$.
Now,  setting
$$
W(z)=\sum_{n=1}^\infty Q_n z^n,
$$
we have  from (\ref{recur})
$$
W(z)= z  \sum_{k \ge 3}w_k(1+W(z))^{2k-2}.
$$
Finally, recalling definition (\ref{wk}) and setting
$$
\phi_E(x)=\sum_{k \ge 3}w_k(1+x)^{2k-2}= \frac{1}{\gamma} \times \frac{ \left( 1- e^{\frac{-1}{\gamma}}\right)^{3}(x+1)^4}{1- \left( 1- e^{\frac{-1}{\gamma}}\right)^{2}(x+1)^2},
$$
we have
$$
W(z)= z \phi_E(W(z)).
$$
We can now use  a well known result in analytic combinatorics (see e.g. \cite{FS},  Proposition IV.5 pag. 278) to conclude  that
$$
Q_n\le \r_\g^n
$$
where
$$
\r_\g=\min_{x>0}{\phi_E(x)\over x}.
$$
An easy computation shows that
for $\gamma=1.569$ or grater we have that $\r_\g < 1$. Therefore  if  $\g\ge 1.569$, recalling (\ref{penne}),
we have that $P_n\le (n+1)^m\r^n_\g$. In other words, the probability that the \textsc{Color-Algorithm} runs  $n$ steps decays exponentially with $n$ as soon as $n$ is sufficiently large and thus the algorithm stops. Thus, the graph $G$ has an acyclic coloring if $N\ge  3.569(\D -1)$.

\vskip.3cm
\\{\bf Acknowledgement:} all authors are partially supported by the brazilian science foundations Conselho Nacional de Desenvolvimento Cient\'\i fico e Tecnol\'ogico (CNPq)  and Coordenação de Aperfeiçoamento de Pessoal de N\'\i vel Superior(CAPES).

\vskip.3cm
\\{\bf Declarations of interest:} none.

\end{document}